\documentclass[titlepage,11pt]{article}
\usepackage{graphicx} 

\usepackage{amsmath,amsthm,amssymb}
\usepackage[margin=1in]{geometry}
\usepackage{booktabs}
\usepackage{multirow}
\usepackage{subcaption}
\usepackage[maxbibnames=5]{biblatex}

\usepackage{authblk}

\newtheorem{proposition}{Proposition}

\title{A Specialized Simplex Algorithm for Budget-Constrained Total Variation-Regularized Problems}
\author[1]{Dominic Yang}
\affil[1]{Argonne National Laboratory, Mathematics and Computer Science Division, dominic.yang@anl.gov}
\begin{document}

\maketitle

\begin{abstract}
    We consider a class of linear programs on graphs with total variation regularization and a budgetary constraint. For these programs, we give a characterization of basic solutions in terms of rooted spanning forests with orientation on the underlying graph. This leads to an interpretation of the simplex method in terms of simple graph operations on these underlying forests. We exploit this structure to produce an accelerated simplex method and empirically show that such improvements can lead to an order of magnitude improvement in time when compared to state-of-the-art solvers.
\end{abstract}

\section{Introduction}

We consider the following problem:

\begin{equation}\label{eq:TV-Original}
\begin{aligned}
    \min_{x} \qquad &\sum_{v \in V} c_v x_v + \alpha\sum_{uv \in E}|x_u - x_v| \\
    \text{s.t.} \qquad &\sum_{v\in V}h_v x_v \le \Delta,\\
                        &x\in [0,1]^{V}
\end{aligned}
\end{equation}
In this setting, we assume that we have an underlying simple directed graph $G = (V,E)$ where $V$ is the set of vertices and $E$ the set of edges. We also assume we have vectors $c, h \in \mathbb{R}^V$. The sole constraint in this problem we refer to as the budget constraint and the corresponding vector $h$ as the representative costs for using resource $x_v$ limited by the upper bound $\Delta$. We generally assume $h_v > 0$ but the results we discuss in this paper can be adapted to $h_v \le 0$ with care.

The second term in the objective is a discrete analog of the Total Variation (TV) of $x$. The TV was first introduced as a regularizer in \cite{rudin1992nonlinear} for the purposes of image denoising where produced solutions generally had large smooth regions with well-defined edges. Recently, problems of this form have emerged in the context of topology optimization \cite{manns2024discrete} where the total variation term has the impact of reducing speckle-patterns in the produced solutions which correspond to easier to construct materials in the physical space. The problem \eqref{eq:TV-Original} therefore may appear as a subproblem where a more complicated nonlinear model is linearized about a given point, and the budgetary constraint exists to ensure that the solution to the subproblem remains close to the point.

To handle \eqref{eq:TV-Original}, we introduce the linearization $|x_u - x_v| = a_{uv} + a_{vu}$, $x_u - x_v = a_{uv} - a_{vu}$ where $a_{uv}$ and $a_{vu}$ are nonnegative variables which represent the amount $x_u$ exceeds $x_v$ and the amount $x_v$ exceeds $x_u$, respectively. In this work, we will actually consider a slightly more general version of the problem where we allow arbitrary coefficients on $a_{uv}$ and $a_{vu}$ in the objective so long as $d_{uv} + d_{vu} \ge 0$.
After introducing slack variable $s$ for the budget constraint, we have the following formulation:
\begin{equation}\label{eq:TV-Linear}
\begin{aligned}
    \min_{x} \qquad &\sum_{v \in V} c_v x_v + \sum_{uv \in E}d_{uv}a_{uv} + d_{vu} a_{vu} \\
    \text{s.t.} \qquad & x_u - x_v = a_{uv} - a_{vu},  \quad uv \in E\\
                        &\sum_{v\in V}h_v x_v + s = \Delta,\\
                        & 0 \le x_v \le 1,\quad  a_{uv}, a_{vu}, s \ge 0, \quad v \in V, uv \in E
\end{aligned}
\end{equation}
We observe that the constraint matrix for this problem is given by
\begin{equation}
A_{TV-B} = \begin{bmatrix}
    A_G^T & -I & I & 0 \\
    h^T & 0 & 0 & 1
\end{bmatrix}
\end{equation}
where $A_G$ is the oriented incidence matrix for $G$. We refer to the matrix $\begin{bmatrix}A_G^T & -I & I \end{bmatrix}$ as $A_{TV}$ noting that it is totally unimodular. We refer to the full constraint matrix for this linear program as $A_{TV-B}$ and the associated polyhedron as $\mathcal{P}_{A_{TV-B}}$.
The constraint matrix's resemblance to a network matrix suggests we may be able to exploit this graphical structure with appropriate care for the budget constraint to develop algorithms for this problem. We show that this intuition is correct and use this structure to characterize the simplex method and accelerate its performance.

\paragraph{Related Work.}
Problems similar to \eqref{eq:TV-Original} without the budget constraint have been known to be solvable efficiently even with integer restrictions. In particular, when the objective is generalized to $\sum_v F_v(x_v) + \sum_{uv} G_{uv}(x_u - x_v)$ and $x$ is restricted to be integral, it is still efficiently solvable, and \cite{ahuja2004cut} and \cite{hochbaum2001efficient} both present minimum cut algorithms for solving the problem. \cite{boykov2002fast} and \cite{boykov2004experimental} address the scenario where $G_{uv}$ is instead an arbitrary function of the two vertices and provide approximation guarantees also by using minimum cut methods. A few papers \cite{zalesky2002network,chambolle2005total} also specifically address the scenario where $G$ is the total variation of a graph as we define it.

In comparison, the case with a budgetary constraint has seen much less attention. Of note are two recent papers by the same authors \cite{severitt2023efficient} and \cite{manns2024discrete} which study the integer version of \eqref{eq:TV-Original} with a 1-norm type constraint instead of the general linear constraint. We will see in this paper that similar results hold for \eqref{eq:TV-Linear} as for the problem they studied.

A problem which bears resemblance to this problem is the constrained maximum flow problem which takes the classical maximum flow problem and adds an additional linear constraint applied the vector of flows. We mention this problem as it is formed by adding a budget constraint to a linear program whose constraint matrix is the oriented incidence matrix, and so is in some sense dual to our problem. First studied in \cite{ahuja1995capacity}, various algorithms have been proposed since then including specialized network simplex algorithms \cite{ccalicskan2011specialized, holzhauser2017network}, capacity scaling algorithms \cite{ccalicskan2009capacity}, and a cost scaling algorithm \cite{ccalicskan2012faster}. See \cite{holzhauser2016budget} and \cite{holzhauser2017complexity} for a comprehensive study of hardness and approximability results.

\paragraph{Terminology and Notation.} 
As we will be working with graphs regularly in this paper, we introduce some terminology and notation to facilitate discussion. Given a subset of vertices $S \subset V$, we denote by $\delta^+(S)$, the set of outgoing edges from $S$ and by $\delta^-(S)$, the set of incoming edges to $S$. The union of these sets we denote by $\delta(S)$. 

Given a vertex $v \in V$, we denote by $e_v$, the vector with $|V|$ components where all components are 0 except for at $v$ which takes value 1. Analogously, for $uv \in E$, $e_{uv}$ is the vector with $|E|$ components where all components are 0 except for at ${uv}$ which takes value 1. Given a subset $S \subset V$ and $T \subset E$, we denote by $\chi_S$ and $\chi_T$, the vectors with $|V|$ and $|E|$ components which are 0 everywhere except for on $S$ and $T$, respectively.
Given a vector $c \in \mathbb{R}^V$ and $d \in \mathbb{R}^E$, we write $c(S) := \sum_{v\in S} c_v$ and $d(T) := \sum_{uv \in T} d_{uv}$.

\paragraph{Outline of Paper.}
In this paper, we elaborate on the structure of basic solutions to \eqref{eq:TV-Linear} and demonstrate how this structure enables us to provide an accelerated simplex algorithm for this class of problems. In Section \ref{sec:polyhedron}, we establish an interpretation of all basic solutions in terms of rooted spanning forests and discuss the various implications this has for the polyhedral structure of \eqref{eq:TV-Linear}. Then, in Section \ref{sec:basic-equations}, we establish equations for all basic variables in terms of nonbasic variables leading to concise formulations of reduced costs and optimality conditions for solutions. In Section \ref{sec:pivoting}, we discuss in detail how the simplex method can be interpreted in terms of the underlying spanning forest representation and use representation to significantly reduce work required. Finally in Section \ref{sec:experiments}, we give a brief experimental study demonstrating significant gains using this accelerated algorithm when compared to the state-of-the-art solver CPLEX.

\section{Polyhedral Structure of $\mathcal{P}_{A_{TV-B}}$}\label{sec:polyhedron}

To begin, we first characterize basic solutions for \eqref{eq:TV-Linear} in terms of rooted spanning forests of $G$.

\begin{proposition} \label{prop:tv-linear-basic}
    Let $x$ be a basic solution of \eqref{eq:TV-Linear}. Let $V_{NB} \subset V$ denote the vertices $v$ where $x_v$ are nonbasic, $E^+$ denote edges $uv \in E$ where $a_{uv}$ is nonbasic, $E^-$ denote edges $uv$ where $a_{vu}$ is nonbasic, $E_{NB} = E^+ \cap E^-$ and $G_{NB} = (V, E_{NB})$ denote the subgraph given by selecting only edges in $E_{NB}$. Let $V= V_1 \cup\dots \cup V_k$ partition the vertices of $G_{NB}$ into connected components. Then the following statements are true:
    \begin{enumerate}
        \item $E^+ \cup E^- = E$.
        \item If $s$ is basic, $|V_{NB} \cap V_i| = 1$ for $i=1,\ldots,k$; if $s$ is nonbasic, $|V_{NB} \cap V_i| = 1$ for all but one $i$ where $|V_{NB} \cap V_i| = 0$ and $h(S_i) \ne 0$.
        \item Each connected component of $G_{NB}$ is a tree.
    \end{enumerate}
\end{proposition}
\begin{proof}
    The first statement follows from the fact that the columns for $a_{uv}$ and $a_{vu}$ are $(-e_{uv}, 0)$ and $(e_{uv}, 0)$ and are dependent so that both can never be in the basis simultaneously.
    For the second statement, we consider two cases. First assume $s$ is basic. Assume that there is some component which contains no nonbasic vertices, and let $S$ be the set of vertices in this component. Let $N(S)$ denote the set of vertices which are not in $S$ but have a neighbor in $S$. As each column corresponding to a vertex variable $v$ is given by
    \[
        A_{TV,v} = (-\sum_{uv \in \delta^-(v)}e_{uv} + \sum_{vw\in \delta^+(v)}e_{vw}, 1),
    \]
    the sum over all such columns is
    \[
        \sum_{v\in S} A_{TV,v} = (-\sum_{wt \in \delta^-(S)}e_{wt} + \sum_{wt \in \delta^+(S)}e_{wt}, h(S))
    \]
    as for each edge $uv$ between vertices in $S$, the corresponding components cancel out. On the other hand, by assumption, one of the variables $a_{uv}, a_{vu}$ for each edge $uv \in \delta(S)$ is basic, and these columns span $\{(e_{uv},1): uv \in \delta(S)\}$. Therefore, we can find a linear combination producing solely the vector $(0, h(S))$ which is linearly dependent with $(0,1)$, the basis column corresponding to $s$. Hence, we do not have a basis, and each component requires at least one nonbasic vertex.

    When $s$ is nonbasic, the column $(0,1)$ is no longer in the basis. If $h(S) = 0$, the argument above still holds and the columns are not linearly independent, but if $h(S) \ne 0$ the above argument does not work as we cannot eliminate the last component. However, if we have two components $S, S'$ both satisfying $h(S), h(S')\ne 0 $, as above we can produce two linear combinations giving $(0, h(S)$ and $(0, h(S'))$ which are linearly dependent. Hence, when $s$ is nonbasic, all but one component must have at least one nonbasic vertex.

    To prove the third statement, we proceed with a counting argument. Note that as there are $|E| + 1$ equations in \eqref{eq:TV-Linear}, and $|V| + 2|E| + 1$ variables, so that there are $|V| + |E|$ nonbasic variables. As $E^+ \cup E^- = E$, we use up $|E|$ nonbasic variables by picking one of $a_{uv}$ and $a_{vu}$ to be nonbasic for each edge $uv \in E$. This leaves $|V|$ variables which correspond to vertices in $V_{NB}$, edges in $E_{NB}$ and possibly $s$. When $s$ is basic, then $|V_{NB}| + |E_{NB}| = |V|$. If $G_{NB}$ has $k$ components, necessarily  $|E_{NB}| \ge |V|-k$, and each component has a nonbasic vertex so that $|V_{NB}| \ge k.$ This is only possible if $|E_{NB}| = |V|-k$ and $|V_{NB}| = k$ implying each component is a tree, and each tree has one nonbasic vertex. When $s$ is nonbasic, then $|V_{NB}| + |E_{NB}| = |V|-1$ and one tree will not have a nonbasic vertex. 
\end{proof}

\begin{figure}
    \centering
    \includegraphics[width=\linewidth]{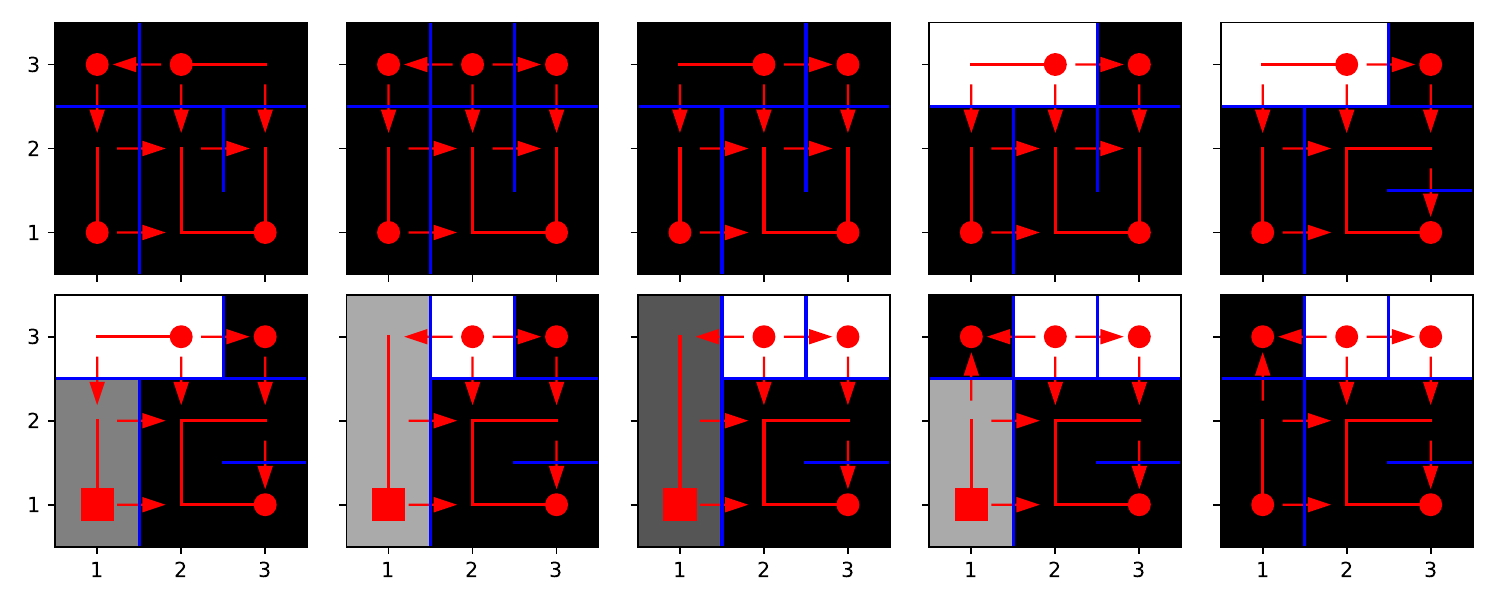}
    \caption{A graph theoretic interpretation of several pivots applied to the linear program \eqref{eq:TV-Linear} where the underlying graph is a $3\times 3$ grid graph, $h_v \equiv 1$, and $\Delta = 3$. Vertices $v$ with a circle are the roots and indicate nonbasic $x_v$ and other vertices correspond to basic $x_v$. Vertices $v$ with a square indicate the root of the basic component which is also basic. Undirected edges are those where both $a_{uv}$ and $a_{vu}$ are nonbasic whereas edges $uv$ with a direction indicate $a_{uv}$ is basic and $a_{vu}$ is nonbasic. The underlying color of a vertex indicates the value of $x_v$ where black signifies $x_v = 0$, white signifies $x_v = 1$, and gray colors indicate intermediate values.}
    \label{fig:budget-simplex-pivots}
\end{figure}

We can also easily characterize which basic solutions are feasible.
\begin{proposition}\label{prop:tv-linear-basic-feasible}
    A basic solution $(x, a, s)$ with $x \in [0,1]^V$ is feasible if $\sum_v h_v x_v \le \Delta$, and for each edge between $u$ and $v$ in different components, if $x_{u} > x_{v}$, $a_{uv}$ is basic.
\end{proposition}
\begin{proof}
    We have $x_u - x_v = a_{uv} - a_{vu} > 0$. As $a_{uv}$ and $a_{vu}$ need to be nonnegative to be feasible and only one may be basic, it must be the case that $a_{uv}$ is basic.
\end{proof}

Propositions \ref{prop:tv-linear-basic} and \ref{prop:tv-linear-basic-feasible} tell us that we can associate any basic solution with a rooted spanning forest of $G$ alongside an orientation of any edge $uv$ between components designating which of $a_{uv}$ or $a_{vu}$ is basic (i.e., if the edge is oriented $uv$, $a_{uv}$ is basic and $a_{vu}$ is nonbasic). For it to be feasible, the orientation needs to be pointed from larger to smaller element. When $s$ is basic, each component has a nonbasic vertex which we designate as the root of the component, and when $s$ is nonbasic one component has no such vertex, and we refer to this component as the \textbf{basic component}. We will still designate a root $r'$ of the basic component, but this can be any vertex in the component. See Figure \ref{fig:budget-simplex-pivots} for depictions of the corresponding rooted spanning forest for several basic solutions. Note that for adjacent vertices $u$ and $v$ in the same component, as both $a_{uv}$ and $a_{vu}$ are nonbasic and thus zero, $x_u = x_v$. It immediately follows that $x$ is constant on each tree in the forest.

We provide some notation to better discuss the concepts associated with these forests. For any given vertex $v$, we let $r(v)$ denote the root of the tree containing $v$. We let $T_v = (N_v, E_v)$ denote the tree comprising of $v$ and all descendants of $v$ with the convention that all edges in $E_v$ are directed away from $v$. We designate vertices $v$ as basic if $x_v$ is basic and nonbasic otherwise. We also designate edges $uv \in E$ as nonbasic if both $a_{uv}$ and $a_{vu}$ are nonbasic, and basic if only one of the variables is nonbasic. 
As $G_{NB}$ is a forest, between two vertices $u,v$ in the same component, there is a unique path from $u$ to $v$ which we denote $P(u,v)$. In this context, we will represent a path by the sequence of edges from $u$ to $v$, i.e., $P(u,v) = (uv_1, v_1v_2, \ldots, v_mv)$. 
Given $x$ corresponding to a basic feasible solution, we will define a \textbf{rooted spanning forest with orientation aligned with $x$} to be any rooted spanning forest with orientation on non-tree edges satisfying $x$ is constant on each tree, and for each edge between trees $uv$, if $x_u > x_v$, then the edge is oriented as $uv$. 

Proposition \ref{prop:tv-linear-basic} leads to the following immediate corollary regarding the polyhedral structure of \eqref{eq:TV-Linear}. A similar result is shown in \cite{manns2024discrete} where instead of our budget constraint, a constraint limiting the distance of $x$ in the 1-norm to a given point is employed.

\begin{proposition}\label{prop:TV-B-Vertices}
    Let $(x, a, s)$ denote a vertex of the underlying polyhedron $\mathcal{P}_{A_{TV-B}}$. Then $x \in \{0,1\}$ everywhere except possibly on a subset $S \subset V$ which induces a connected subgraph of $G$.
\end{proposition}
\begin{proof}
    From Proposition \ref{prop:tv-linear-basic}, we can associate $(x, a, s)$ with a (not necessarily unique) spanning forest where all but at most one tree has a designated root. For each rooted tree, the value $x_v$ for the nonbasic vertex dictates the value for the entire tree and as $x_v$ is nonbasic it takes value in $\{0, 1\}$. Only on vertices of the basic component which is connected can $x$ take a different value. 
\end{proof}

Another significant consequence of Proposition \ref{prop:tv-linear-basic} is that there are an extreme number of basic solutions in comparison to the number of vertices. In particular, for every $x \in \{0,1\}^V$, every single choice of rooted spanning forest with orientation aligning with $x$ corresponds to a basic solution. We can compute the number of bases using the following proposition.
\begin{proposition}\label{prop:number-of-bases}
    Given a vertex $(x, a, s)$ of $\mathcal{P}_{A_{TV-B}}$, let $G_x = (V, E_x)$ be the subgraph of $G$ where $uv \in E$ is in $E_x$ if $x_u = x_v$, and let $A_{G_x}$ be the oriented incidence matrix of $G_x$. The number of basic feasible solutions corresponding to $x$ is between $|\det(A_{G_x}^TA_{G_x} + 2I)|$ and $2|\det(A_{G_x}^TA_{G_x} + 2I)|$.
\end{proposition}
\begin{proof}
    Since $A_{TV}$ is a totally unimodular, by Theorem 3 in \cite{maurer1976matrix}, we can compute the number of bases of $A_{TV}$ as $|\det({A_{TV}A_{TV}^T})| = |\det(A_{G}^TA_G + 2I)|$, and these bases all correspond to rooted spanning forests with orientation on $G$. To ensure that there are no tree edges between vertices $u$ and $v$ where $x_u$ and $x_v$ take different value, we restrict the graph to $G_x$ which ensures that the spanning forests align with $x$. The number of rooted spanning forests where one component lacks a root is clearly upper bounded by the number of rooted spanning forests. Hence, the total number of basic solutions corresponding to $(x, a, s)$ is between $|\det(A_{G_x}^TA_{G_x} + 2I)|$ and $2|\det(A_{G_x}^TA_{G_x} + 2I)|$.
\end{proof}

As an example of the worst case of how fast this quantity grows, when $G_x$ is the complete graph on $n$ nodes, $|\det(A_{G_x}^TA_{G_x} + 2I)| = 2(n-2)^{n-1}$. To see this note, $L = A_{G_x}^TA_{G_x}$ is the Laplacian of the complete graph, and $L = B - nI$ where $B$ is the $n\times n$ matrix of all ones. Observe the vector of all ones is an eigenvector of $L + 2I$ with eigenvalue 2, and $e_u - e_v$ for vertices $u$ and $v$ is an eigenvector with eigenvalue $2-n$ of which we can find $n-1$ linearly independent eigenvectors. Proposition \ref{prop:number-of-bases} indicates that the underlying linear program is extremely degenerate, and we should expect a large number of degenerate pivots should we attempt to solve \eqref{eq:TV-Linear} using the simplex method. We will see in Section \ref{sec:experiments} that the simplex method performs well in spite of this fact.

\section{Linear Program Structure}\label{sec:basic-equations}
In this section, we examine further the structure of the linear program and derive equations for each of the basic variables in terms of the nonbasic variables. From these equations, we derive succinct formulations for the reduced costs of nonbasic variables and establish optimality conditions for \eqref{eq:TV-Linear}.

\subsection{Basic Variable Equations for Basic $s$}

We will first consider the case where $s$ is a basic vertex.
Given the spanning forest structure of the basic solutions, we can devise simple formulas for the basic variables in terms of nonbasic variables. 
Since we have an explicit formulation for $x_u - x_v$, we write basic $x_v$ in terms of the nonbasic root vertex $x_{r(v)}$ and the path of nonbasic edge variables between $v$ and $r(v)$. In particular, we have
\begin{equation}\label{eq:basic-x}
    x_v = x_{r(v)} - \sum_{wt \in P(r(v),v)}(x_w - x_t) = x_{r(v)} - \sum_{wt \in P(r(v),v)}(a_{wt} - a_{tw}).
\end{equation}
For basic edge variables, $a_{uv} = a_{vu} + x_u - x_v$, where $a_{vu}$ is nonbasic, we can substitute $x_u$ and $x_v$ using \eqref{eq:basic-x} to get a formulation purely in terms of nonbasic variables. We have
\begin{align}
    a_{uv} = a_{vu} + x_{r(u)} - \sum_{wt \in P(r(u),u)}(a_{wt} - a_{tw}) - x_{r(v)} + \sum_{wt \in P(r(v),v)}(a_{wt} - a_{tw}), \label{eq:basic-a-plus}\\
    a_{vu} = a_{uv} - x_{r(u)} + \sum_{wt \in P(r(u),u)}(a_{wt} - a_{tw}) + x_{r(v)} - \sum_{wt \in P(r(v),v)}(a_{wt} - a_{tw}).\label{eq:basic-a-minus}
\end{align}
In the event that $r(u) = r(v)$, we can simplify the above expression by canceling out any terms corresponding to edges which appear in both $P(u,r(u))$ and $P(v, r(v))$. This lets us write the expression as a sum over the path from $u$ to $v$, $P(u, v)$:
\begin{equation}\label{eq:basic-a-plus-loop}
    a_{uv} = a_{vu} + \sum_{wt \in P(u,v)}(a_{wt} - a_{tw}), \quad    a_{vu} = a_{uv} - \sum_{wt \in P(u,v)}(a_{wt} - a_{tw}). 
\end{equation}
Finally for $s$, we have the following formulation given by substituting $x$ using \eqref{eq:basic-x}:
\begin{equation}\label{eq:basic-s}
    s = \Delta - \sum_{v \in V}h_vx_v = \Delta - \sum_{r \in V_{NB}}\left(h(N_r)x_r - \sum_{uv \in E_r}h(N_v)(a_{uv} - a_{vu}) \right).
\end{equation}

\paragraph{Reduced Costs}
Equations \eqref{eq:basic-x} through \eqref{eq:basic-s} immediately tell us several things about how the simplex algorithm operates for the linear program \eqref{eq:TV-Linear}. For each nonbasic root vertex $r$, we observe from \eqref{eq:basic-x} that attempting to raise the value of $x_r$ corresponds to an equal increase in value for each other basic $x_v$ where $r$ is the root for $v$. We also see that each basic $a_{uv}$ where $r(u) = r$ increases in value and for each basic $a_{uv}$ where $r(v) = r$ decreases in value, i.e., for all edges leaving $N_r$, the basic variable increases, and for edges entering $N_r$, the basic variable decreases. 
Hence, increasing $x_r$ corresponds to uniformly increasing the value of $x$ on $N_r$ as well as $a$ on $\delta^+(N_r)$ and a uniform decrease of $a$ on $\delta^-(N_r)$.
Altogether, the reduced cost for $x_r$ is given by
\begin{equation}\label{eq:reduced-cost-vertex}
    r_{x_r} = \sum_{v \in N_r} c_v + \sum_{uv \in \delta^+(N_r)}d_{uv} - \sum_{uv \in \delta^-(N_r)}d_{uv} = c(N_r) + d(\delta^+(N_r)) - d(\delta^-(N_r)).
\end{equation}
As the expression $c(S) + d(\delta^+(S) - d(\delta^-(S))$ will show up repeatedly, we denote this quantity by $F(S)$ for $S \subset V$.

We can make a similar argument to determine reduced costs for the variables $a_{uv}$ for a given edge $uv$. When $uv$ is in a tree $E_r$ (meaning both $a_{uv}$ and $a_{vu}$ are nonbasic and $uv$ is directed away from $r$), $a_{uv}$ appears in equations for basic variables $x_{w}$ and $a_{wt}$ whenever $uv$ is in the paths $P(r(w), w)$ or $P(r(t), t)$. The set of all such vertices $w$ for which this is true is $v$ and all descendants of $v$, $N_v$, and from $\eqref{eq:basic-x}$, we see that raising the value of $a_{uv}$ corresponds to an equal decrease in $x_w$. Similarly, the set of edges $wt$ where $uv$ is in the path to $w$ are the edges leaving $N_v$, $\delta^+(N_v)$, and the edges where $uv$ is in the path to $t$ are those edges entering $N_v$, $\delta^-(N_v)$. In the former case, we observe that increasing $a_{uv}$ decreases $a_{wt}$ as in \eqref{eq:basic-a-plus}, and in the latter case the opposite occurs as in \eqref{eq:basic-a-minus}. Altogether, this gives the following formula for the reduced cost of $a_{uv}$:
\begin{equation}\label{eq:reduced-cost-edge-plus}
    r_{a_{uv}} = d_{uv} - \sum_{v \in N_v}c_v - \sum_{wt \in \delta^+(N_v)}d_{wt} + \sum_{wt \in \delta^-(N_v)}d_{wt} = d_{uv} - F(N_v). 
\end{equation}
Analogously, for $a_{vu}$, since in \eqref{eq:basic-x}, \eqref{eq:basic-a-plus}, \eqref{eq:basic-a-minus}, the sign is flipped, we have the following expression:
\begin{equation}\label{eq:reduced-cost-edge-minus}
    r_{a_{vu}} = d_{vu} + F(N_v). 
\end{equation}
For the case of the edge $vu$ where $a_{uv}$ is basic, we see that $a_{uv}$ only appears in \eqref{eq:basic-a-plus} with coefficient 1 giving the following expression for the reduced cost:
\begin{equation}\label{eq:reduced-cost-basic-edge}
    r_{a_{uv}} = d_{uv} + d_{vu}.
\end{equation}
Note that by assumption, we have $d_{uv} + d_{vu} \ge 0$, and therefore, we never pivot on these variables.

\subsection{Basic Variable Equations for Nonbasic $s$}

When $s$ is nonbasic, one of the trees in the spanning forest representation will not have a nonbasic vertex. Recall we still designate one vertex as root $r'$. From \eqref{eq:basic-s}, we can first derive an expression for $x_{r'}$ as follows (recall, the tree without a nonbasic vertex satisfies $h(N_{r'}) \ne 0$):
\begin{equation}
    h(N_{r'})x_{r'} = \Delta - s - \sum_{r \text{ root}}\left(h(N_r)x_r - \sum_{uv \in E_r}h(N_v)(a_{uv} - a_{vu})\right)  + \sum_{uv \in E_{r'}} h(N_v)(a_{uv} - a_{vu})
\end{equation}

If we substitute this into the expressions \eqref{eq:basic-x} for $w \in N_{r'}$, we get the following:

\begin{align}
    x_w    &= \frac{\Delta - s}{h(N_{r'})} -\sum_{r \text{ root}}\left(\frac{h(N_r)}{h(N_{r'})}x_r - \sum_{uv \in E_r}\frac{h(N_v)}{h(N_{r'})}(a_{uv} - a_{vu})\right)  \label{eq:basic-x-s}\\
        &\quad+ \sum_{uv \in E_{r'}\setminus P(r',w)} \frac{h(N_v)}{h(N_{r'})}(a_{uv} - a_{vu}) + \sum_{uv \in P(r',w)}\frac{h(N_v)-h(N_{r'})}{h(N_{r'})}(a_{uv} - a_{vu}) \nonumber
\end{align}

For basic edges $wt$ where $w$ is in the basic component, we have the following expression for $a_{wt}$ (and $a_{tw}$ can be determined analogously):
\begin{equation} \label{eq:basic-a-plus-s}
    a_{wt} = a_{tw} + x_{r'} - \sum_{uv \in P(r', w)}(a_{uv}-a_{vu}) - x_{r(t)} + \sum_{uv \in P(r(t), t)}(a_{uv} - a_{vu})
\end{equation}
This equation informs us that $a_{wt}$ for $wt \in \delta^+(N_{r'})$ rises at the same rate as $x_{r'}$ and $a_{wt}$ for $wt \in \delta^-(N_{r'})$ decreases at that rate.
All other variable equations remain as before.

\paragraph{Reduced Costs.}


Again, we compute the reduced costs for each nonbasic variable. Importantly, we note that as before, shifting $x_r$ and $a_{uv}$ where $r$ and $uv$ are in the rooted components corresponding to shifting $x$ on $N_r$ and $N_v$, respectively. However, from equations \eqref{eq:basic-x-s} and \eqref{eq:basic-a-plus-s}, we see that such shifts also impact the basic component as well. 
From \eqref{eq:basic-x-s}, we see that raising $x_r$ corresponds to decreasing $x_w$ at a rate of $h(N_r)/h(N_{r'})$ for each $w$ in the basic component. Likewise for each each $wt \in \delta^+(N_{r'})$, we decrease $a_{wt}$, and for each $wt \in \delta^-(N_{r'})$, we increase $a_{wt}$, both at a rate of $h(N_{r})/h(N_{r'})$. Altogether, we then have
\begin{equation}\label{eq:reduced-cost-vert-s}
    r_{x_r} = F(N_r) + \frac{h(N_r)}{h(N_{r'})} \left(-\sum_{v \in N_{r'}}c_v - \sum_{uv \in \delta^+(N_{r'})}d_{uv} + \sum_{uv \in \delta^-(N_{r'})}d_{uv}\right) = F(N_r) - \frac{h(N_r)}{h(N_{r'})}F(N_{r'}).
\end{equation}


For nonbasic edges $uv \in E_r$, similarly, we maintain that same shift on $N_v$ as before, and we also observe a change in the opposite direction on $N_{r'}$ with a rate of change of $h(N_v)/h(N_{r')}$. Therefore,
\begin{equation}\label{eq:reduced-cost-edge-plus-s}
    r_{a_{uv}} = d_{uv} - F(N_v) + \frac{h(N_v)}{h(N_{r'})}F(N_{r'}), \quad r_{a_{vu}} = d_{vu} + F(N_v) - \frac{h(N_v)}{h(N_{r'})}F(N_{r'}).
\end{equation}

For edges $uv$ in the basic component, observe that if $u$ is on the same side as $r'$, from \eqref{eq:basic-x-s}, raising $a_{uv}$ raises $x_w$ by $(h(N_v) - h(N_{r'}))/h(N_{r'}) = h(N_{r'}\setminus N_v)/h(N_{r'})$ for $w \in N_v$, and increases by the same amount for $a_{wt}$ for $wt \in \delta^+(N_v)$ and decreases for $wt \in \delta^-(N_v)$. All other vertices and edges in $N_{r'}\setminus N_v$ and $\delta(N_{r'}\setminus N_v)$ change at a rate of $h(N_v)/h(N_{r'})$. We then have
\begin{equation}
    r_{a_{uv}} = d_{uv} - \frac{h(N_{r'}\setminus N_v)}{h(N_{r'})}F(N_v) + \frac{h(N_v)}{h(N_{r'})}F(N_{r'}\setminus N_v).
\end{equation}
For these edges, we see that we would be incurring the cost for decreasing $N_v$ and increasing $N_{r'} \setminus N_v$ simultaneously if we pivot on $a_{uv}$ and vice versa for $a_{vu}$.

Lastly, we see that $s$ appears in the \eqref{eq:basic-x-s} and \eqref{eq:basic-a-plus-s} for all $v \in N_{r'}$ and $uv \in \delta(N_{r'})$ with a coefficient of $-1/h(N_{r'})$. Hence, we have the reduced cost is given by
\begin{equation}\label{eq:reduced-cost-s}
    r_s = \frac{1}{h(N_{r'})}\left(-\sum_{v \in N_{r'}}c_v - \sum_{uv \in \delta^+(N_{r'})}d_{uv} + \sum_{uv \in \delta^-(N_{r'})}d_{uv}\right) = -\frac{F(N_{r'})}{h(N_{r'})}.
\end{equation}
Hence, pivoting $s$ back into the basis corresponds to reducing the component $N_{r'}$ back to 0.

\subsection{Optimality Conditions}
Knowing the equations for reduced costs for all the nonbasic variables allows us to succinctly state optimality conditions for \eqref{eq:TV-Linear}. Given that we have a minimization problem, optimality occurs exactly when all the reduced costs for the nonbasic variables at their lower bound (i.e., 0) are nonnegative and the reduced costs for nonbasic variables at their upper bound (when $x_r = 1$) are nonpositive. We can then establish characterizations of optimality in terms of these rooted spanning forests. As before, we divide these conditions into two cases, when $s > 0$ and when $s = 0$.

\begin{proposition}\label{prop:optimality-sub-budget}
    A feasible point $(x, a, s)$ for \eqref{eq:TV-Linear} with $\sum_{v}h_vx_v < \Delta$ is optimal if $x \in \{0,1\}^V$ and there exists a rooted spanning forest with orientation aligned with $x$ such that $F(x_r) \ge 0$ for roots $r$ where $x_r = 0$, $F(x_r) \le 0$ for roots $r$ with $x_r = 1$, and for each tree edge $uv$, $d_{uv} - F(N_v) \ge 0$ and $d_{vu} + F(N_v) \ge 0$.
\end{proposition}
In other words, Proposition \ref{prop:optimality-sub-budget} states if we can find a rooted spanning forest with orientation of non-tree edges where we cannot shift $x$ on any tree or subtree to produce a better solution, then we are optimal.  For the case where $s = 0$, as $h(S) > 0$ for all $S \subset V$, we can consider instead the weighted reduced costs $r_{x_r}/h(N_r)$ and $r_{a_{uv}}/h(N_v)$. This leads to the following proposition:
\begin{proposition}\label{prop:optimality-equal-budget}
    A feasible point $(x, a, s)$ for \eqref{eq:TV-Linear} with $\sum_v h_vx_v = \Delta$ is optimal if $x \in \{0, 1\}^V$ except possibly one connected component $S$ where $x$ is constant, and there exists a rooted spanning forest with orientation aligned with $x$ such that $F(x_r)/h(x_r) \ge F(S)/h(S)$ for each root $r$ with $x_r = 0$, $F(x_r)/h(x_r) \le F(S)/h(S)$ for each root $r$ with $x_r = 1$,  for each tree edge $uv$, $(F(N_v) - d_{uv})/h(N_v) \le F(S)/h(S)$ and $(d_{vu} + F(N_v))/h(N_v) \ge F(S)/h(S)$, and $F(S)/h(S) \le 0$.
\end{proposition}
Proposition \ref{prop:optimality-equal-budget} can be understood to say that a point $x$ is optimal if we can find a rooted spanning forest with orientation where $F(S) \le 0$ any move which shifts $x$ on a tree or subtree would lead to less efficient usage of the budget constraint than $S$.

\section{Pivoting}\label{sec:pivoting}

In this section, we recount how pivots in the simplex method can be interpreted under the rooted spanning forest representation of basic solutions given in Proposition \ref{prop:tv-linear-basic}. Pivoting can be understood as transitioning between various rooted spanning forests of the underlying graph. 
Introducing a nonbasic edge variable $a_{uv}$ to the basis signifies splitting a tree at a given edge and orienting the edge in the direction $uv$ whereas pivoting a basic edge out of the basis means merging two adjacent trees along a given edge. 
Similarly, pivoting a variable $x_v$ in and out of the basis means marking or unmarking $v$ as a root of a tree.
As discussed in the previous section, all such pivots also attempt to shift the values of $x$ on specific subtrees of the spanning forest.
We catalog the variety of ways in which these merges, splits, and shifts can occur in the simplex method and eventually reach a forest satisfying the conditions of Propositions \ref{prop:optimality-sub-budget} and \ref{prop:optimality-equal-budget}.

Fundamental to our discussion is the notion of a blocking edges. Analogous to the concept in network simplex methods, when we attempt to shift a region $N_v$ up, we define a \textbf{blocking edge} as any incoming edge $uv$ to $N_v$ where $x_u - x_v$ is minimized (analogously if shifting down, it is any such outgoing edge $uv$ from $N_v$). These edges are integral in that if we shift $x$ too far on these vertices (further than $x_u - x_v$), then we produce an infeasible solution as $x$ will not be aligned with the orientation as in Proposition \ref{prop:tv-linear-basic-feasible}. As such, these blocking edges then signify how far we can shift as well as which basic variables $a_{uv}$ we can pivot out of the basis (and therefore merge along). In practice, there will be several such choices of blocking edges, a consequence of the extreme degeneracy in the problem, and any choice is permissible.

\begin{figure}
    \centering
    \begin{subfigure}[t]{0.32\textwidth}
        \includegraphics[]{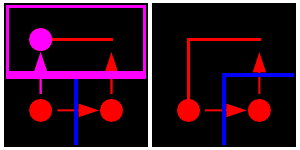}
        \caption{} \label{fig:merge}
    \end{subfigure}
    \begin{subfigure}[t]{0.32\textwidth}
        \includegraphics[]{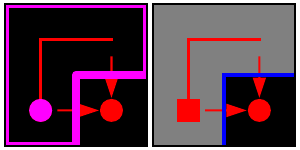}
        \caption{}\label{fig:shift-budget}
    \end{subfigure}
    \begin{subfigure}[t]{0.32\textwidth}
        \includegraphics[]{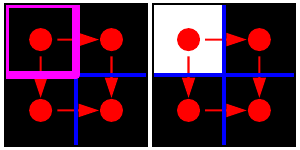}
        \caption{}\label{fig:shift}
    \end{subfigure}
    \begin{subfigure}[t]{0.32\textwidth}
        \includegraphics[]{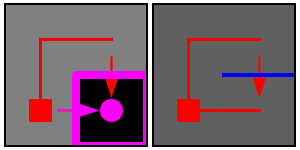}
        \caption{}\label{fig:double-shift-merge}
    \end{subfigure}
    \begin{subfigure}[t]{0.32\textwidth}
        \includegraphics[]{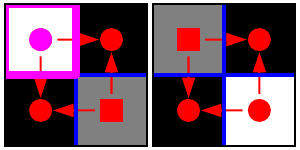}
        \caption{}\label{fig:double-shift}
    \end{subfigure}
    \begin{subfigure}[t]{0.32\textwidth}
        \includegraphics[]{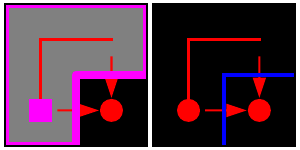}
        \caption{}\label{fig:shift-slack}
    \end{subfigure}
    \begin{subfigure}[t]{0.49\textwidth}
        \includegraphics[]{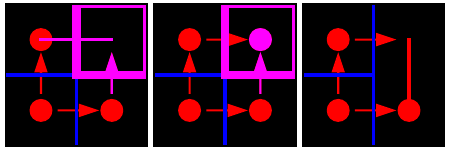}
        \caption{}\label{fig:split}
    \end{subfigure}
    \begin{subfigure}[t]{0.49\textwidth}
        \includegraphics[]{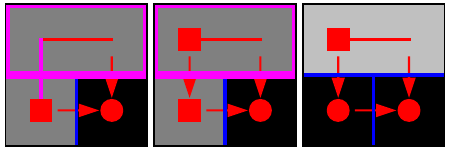}
        \caption{}\label{fig:split-double-shift}
    \end{subfigure}
    \caption{Graphic representation of pivots for \eqref{eq:TV-Linear} where $G$ is a $2\times 2$ grid, $h_v \equiv 1$ and $\Delta = 3/2$. Highlighted circles, undirected edges, and squares indicate entering nonbasic variables, and highlighted directed edges indicate blocking edges. The outlined region indicates vertices to shift.}
    \label{fig:pivots}
\end{figure}

\subsection{Pivoting on Vertices}
As discussed in Section \ref{sec:basic-equations}, pivoting on a vertex variable $x_r$ corresponds to an attempt to shift the values of $x$ on the vertices of the tree $T_r$. We divide our discussion into the cases where $s$ is in the basis and there is no basic component and when $s$ is not in the basis.

\paragraph{Basic $s$.}
When $s$ is in the basis, pivoting in a nonbasic variable $x_r$ attempts to shift the vertices $N_r$. When shifting, there are three scenarios where we can become infeasible: when there is a blocking edge on $N_v$, when shifting will exceed the budget constraint, and when shifting leaves the bounds of $[0,1]$ for $x$. In the first case, simply merge along any blocking edge $uv$ pivoting $x_r$ into the basis and pivoting $a_{uv}$ out of the basis as in Figure \ref{fig:merge}. In the second, we shift until we reach the budget in which case $x_r$ enters the basis and $s$ leaves the basis making it so that $N_r$ is the basic component as in Figure \ref{fig:shift-budget}. In the last case, as we are not concerned with blocking edges or the budget constraint, we can shift up until either $0$ or $1$ and we do not need to perform a basis exchange as can be seen in Figure \ref{fig:shift}.

\paragraph{Nonbasic $s$.}
When $s$ is not in the basis, each pivot on a nonbasic vertex variable $x_r$ shifts $N_r$ and also the basic component $N_{r'}$ in the opposite direction at $h(N_r)/h(N_{r'})$ times the rate. As before, we need to be concerned about blocking edges and the boundaries of $0$ and $1$ for $N_r$, but we also need to ensure that $x$ on the basic component $N_{r'}$ does not leave $[0,1]$. From this, we need to consider to additional cases where we may become infeasible, when there is a blocking edge with $N_{r'}$ and when $x$ on $N_{r'}$ reaches either $0$ or $1$. In the former case with edge $uv$, we shift both components until they are equal, pivoting $x_r$ into the basis and $a_{uv}$ out as can be seen in Figure \ref{fig:double-shift-merge}. In the latter case, we shift until the basic component reaches the boundary and pivot the variable $x_{r'}$ out of the basis to become the root, and we pivot $x_r$ into the basis. This can be seen in Figure \ref{fig:double-shift}. We also mention the case where we pivot $s$ into the basis corresponding to shifting $x$ on $N_{r'}$ down to 0 and introducing $r'$ as the root which can be seen in Figure \ref{fig:shift-slack}.

\subsection{Pivoting on Edges}
When we pivot on a tree edge $a_{uv}$ for $uv$ oriented away from the root, we split on edge $uv$ and are attempting to shift the subtree with vertices $N_v$ down. Analogously, pivoting $a_{vu}$ into the basis also splits the edge and attempts to shift $N_v$ up. In almost all cases, these pivots can be understood as first splitting the tree into two trees, the second with root $v$, and then attempting to shift the new root variable $x_v$ in the appropriate direction as can be seen in Figure \ref{fig:split}. One case which deserves special comment is when we pivot on edge $uv$ in the basic component. In this case, we split along this edge producing two components shifting the component with $v$ down and the component with $u$ up at the rate which keeps the budget constraint tight. Once one component reaches $0$ or $1$, we mark the root as nonbasic for that component and the other component becomes the basic component. This process can be seen in Figure \ref{fig:split-double-shift}.

\subsection{Optimizations}
The simplex method can be significantly accelerated using a few observations. In particular, from Propositions \ref{prop:optimality-sub-budget} and \ref{prop:optimality-equal-budget}, we see that we can easily check for optimality provided we know $y_v := F(N_v)$ and $z_v := h(N_v)$ for each vertex $v$. For these quantities, the following relations hold:
\begin{equation}\label{eq:y-and-z}
    y_v = c_v + d(\delta^+(\{v\})) - d(\delta^-(\{v\})) + \sum_{u \text{ child  of }v}y_u, \quad z_v = h_v + \sum_{\text{ u child of }v}z_u.
\end{equation}
Hence, assuming that the descendants of $v$ are valid, $y_v$ and $z_v$ can be computed by adding up values for direct children. In addition, as these values only depend on descendants, whenever a split occurs at $uv$ creating a new subtree with root $r_v$, only ancestors of $u$ and the ancestors of $v$ in the new tree need to have $y$ and $z$ updated. These are purely those on the path from $r(u)$ to $u$ and the path from $v$ to $r_v$. Similarly, when a merge across edge $uv$ occurs, only vertices from on the path from $r(u)$ to $u$ and from $r(v)$ to $v$ need to have $y$ and $z$ updated. Using equations \eqref{eq:y-and-z}, this can be done in a single pass over these paths.
If $P_{\max}$ is the maximal length of path from root to leaf in any of the trees and $deg_{\max}$ is the maximal degree of any vertex, then such an operation can be done in $O(P_{\max} deg_{\max})$ time. In the worst case, this is $O(|V||E|)$ but in practice, it is significantly faster than simply recomputing all $y_v$ and $z_v$.

The time to perform a full pivot is comprised of the time to find a blocking edge, the time to potentially shift an entire component and the basic component, and finally the time to update the variables $y$ and $z$. Finding a blocking edge requires in the worst case traversing all nodes $N_v$ and edges in $\delta(N_v)$, $O(N_{\max} + \delta_{\max})$ where $N_{\max}$ and $\delta_{\max}$ are maximal values of $|N_v|$ and $|\delta(N_v)|$. Shifting requires in the worst case $O(N_{\max})$ work. Altogether this means a pivot requires $O(P_{\max}deg_{\max} + N_{\max} + \delta_{\max})$ operations in the worst case.

\section{Experiments}\label{sec:experiments}

In this section, we expound on a brief study on the efficiency of the proposed adaptation of the simplex algorithm when compared to the state-of-the-art linear program solver CPLEX. CPLEX provides several different approaches for solving linear programs, and we focus on its simplex implementations, specifically the primal and dual simplex algorithms, alongside its barrier method, an implementation of an interior point algorithm for solving linear programs.

For our study, we randomly generate instance of \eqref{eq:TV-Original} where the underlying graph is a 2-dimensional grid graph, that is, $G = (V, E)$ where $V = \{(i, j) : i, j=1,\ldots,N\}$ for a given $N$ and two vertices $(i, j)$ and $(i', j')$ are adjacent to each other if $|i - i'| + |j - j'| = 1$. In this way, these problems represent total variation-regularized optimization on image processing problems where the vertices represent pixels and edges exist for each pair of adjacent pixels in the image. For the budget constraint, we weight each vertex equally, i.e., we take $h_v \equiv 1$ so that the constraint is simply $\sum_v x_v \le \Delta$. To decide the coefficients $c_v$, we set $c_v = t_v - \epsilon$ where $t_v$ is distributed according to a standard normal distribution and $\epsilon$ is a small constant chosen to ensure that the zero solution is not optimal. 
We choose $N \in \{64, 128, 256, 512\}$ and $\alpha \in \{0.5, 1, 2\}$ and for each choice of $(N, \alpha)$, we construct 10 distinct instances with $c$ as chosen above. To determine $\Delta$, we solve the problem without the budget constraint and set $\Delta$ to be half the budget used by the solution to the unconstrained problem. This way, the budget constraint has an impact on the optimal solution.

Our experiments are performed on an Intel(R) Xeon(R) Gold 6140 CPU with 1621 Gigabytes of memory. We use CPLEX version 20.1.0.0 with all the default settings except that we specify which method we are using (primal simplex, dual simplex, or barrier method), we restrict each run to a single thread, and we set a timeout of one hour (3600 seconds). 

\begin{table}[]
    \centering
\begin{tabular}{llllrrrr}
\toprule
$N$ & $|V|$ & $|E|$ & $\alpha$ & Our Method & Primal Simplex & Dual Simplex & Barrier Method \\
\midrule
\multirow{3}{*}{64} & \multirow{3}{*}{4096} & \multirow{3}{*}{8064} & 0.5 & 0.07 & 0.25 & 0.38 & 0.34 \\
 & & & 1 & 0.14 & 0.45 & 0.85 & 0.41 \\
 & & & 2 & 0.26 & 0.96 & 0.46 & 0.43 \\
\cline{1-8}
\multirow{3}{*}{128} & \multirow{3}{*}{16384} & \multirow{3}{*}{32512} & 0.5 & 0.30 & 1.93 & 4.85 & 1.84 \\
 & & & 1 & 0.75 & 4.85 & 16.32 & 2.38 \\
 & & & 2 & 2.04 & 19.95 & 8.22 & 5.87 \\
\cline{1-8}
\multirow{3}{*}{256} & \multirow{3}{*}{65536} & \multirow{3}{*}{130560} & 0.5 & 1.35 & 24.30 & 97.33 & 10.95 \\
 & & & 1 & 4.29 & 43.72 & 312.84 & 14.66 \\
 & & & 2 & 21.31 & 406.89 & 300.71 & 154.54 \\
\cline{1-8}
\multirow{3}{*}{512} & \multirow{3}{*}{262144} & \multirow{3}{*}{523264} & 0.5 & 5.88 & 340.16 & 1592.66 & 64.74 \\
 & & & 1 & 21.96 & 524.29 & 3607.02 & 85.65 \\
 & & & 2 & 243.86 & 3324.55 & 3607.06 & 2797.20 \\
\cline{1-8}
\bottomrule
\end{tabular}
    \caption{Time (s) to solve a TV-regularized problem with $c_v \sim N(0,1)$ on an $N\times N$ grid graph for given $\alpha$ averaged over 10 choices of $c$.}
    \label{tab:simplex-results}
\end{table}

In Table \ref{tab:simplex-results}, we present the time to determine the optimal solution for each approach for each choice of $N$ and $\alpha$ averaged over ten random choices for $c$. We also include the number of vertices and the number of edges for each graph.

By and large, we observe that our specialized simplex approach for this problem outperforms each other method in CPLEX by a significant margin. In the smaller examples, our approach finishes three to four times as fast but for some of the largest examples, it outperforms the other methods by a factor of 10 to 20. Out of the CPLEX approaches, the barrier method seems to perform the best, and we suspect that this is in large part due to its ability to completely circumvent the degeneracy of the underlying linear program using an interior point approach. In many scenarios, we observed that the primal simplex and dual simplex approaches did not terminate within an hour for the largest instances. As $\alpha$ increases, we observe that solve times become longer across the board, but our approach still maintains a significant edge.

Overall, the results clearly demonstrate that the optimizations that we perform to reduce the amount of work needed in the simplex method has the effect of significantly reducing time to completion in some cases by an order of magnitude when compared to CPLEX.


\section{Conclusion}
In this paper, we develop a characterization of the simplex algorithm for budget-constrained total variation-regularized linear programs. In a manner similar to network simplex algorithms, we characterize basic solutions as corresponding to a rooted spanning forests, and simplex pivots as corresponding to merging, splitting, and shifting trees in this forest. This leads to efficient update rules for the simplex method in terms of these graph operations which allows us to significantly improve the solve time for the simplex method on these problems. Empirically, we observe that when compared with state-of-the-art linear program solvers, our method can outperform these other approaches sometimes by an order of magnitude.

Future research for problems of this type involves a more in-depth analysis on the number of pivots for these kinds of problems. There are strongly polynomial network simplex methods, and as this is an extension of a network-type program, it is possible that such methods could extend to this problem. The extreme degeneracy of the problem also could give rise to randomized methods for solving these problems. Choosing the blocking edge randomly would lead to an algorithm which randomly merges and splits these trees in a manner similar to Karger's algorithm \cite{karger1993global}. Studying the behavior of such randomized algorithms could also provide an interesting line of research.

\paragraph{Acknowledgments}
This work was supported by the U.S. Department of Energy, Office of Science, Advanced Scientific Computing Research, under contract number DE-AC02-06CH11357.

\printbibliography

\end{document}